\documentstyle[amsfonts,12pt]{article}
\newtheorem{theorem}{\sc Theorem}
\newtheorem{proposition}[theorem]{\sc Proposition}
\newtheorem{lemma}[theorem]{\sc Lemma}

\newcommand{\cha}{cardinality $\ha_1$}
\newcommand{\hao}{\ha_1}
\newcommand{\goo}{\go_1}
\newcommand{\restriction}{{|}} 

\newcommand{\proof}{\noindent {\sc Proof. }} 
 
\def\rest{\mathord{\restriction}}

\newcommand{\open}{\Bbb}
 
\def\deq{\mathop=\limits^{\rm def}} 

\newlength{\labparwidth}
\newcommand{\labpar}[2]{$$\parbox{\labparwidth}{#2}\leqno(#1)$$}

\newcommand{\dom}{\mbox{\rm dom}}

\newcommand{\ann}{\mbox{\rm Ann}}

\newcommand{\supp}{\mbox{\rm supp}}
\newcommand{\se}{\subseteq}
\newcommand{\set}[2]{\{#1 \colon #2\}} 
\newcommand{\fin}{$\Box$\par\medskip} 
\newcommand{\qed}{$\Box$\par\medskip}

\newcommand{\ga}{\alpha}
\newcommand{\gb}{\beta}

\newcommand{\gd}{\delta}

\newcommand{\gth}{\theta}

\newcommand{\gk}{\kappa}

\newcommand{\gn}{\nu}

\newcommand{\gr}{\rho}
\newcommand{\gs}{\sigma}
\newcommand{\gt}{\tau}

\newcommand{\go}{\omega}
 
\newcommand{\gG}{\Gamma}

\newcommand{\ha}{\aleph}
 
\newcommand{\oP}{{\open P}}
\newcommand{\oQ}{{\open Q}}

\def\psupp{\mathop{\hbox{p-supp}}}

\title{Explicitly Non-Standard Uniserial Modules}
\author{P. C. Eklof \\University of California, Irvine\thanks{Thanks
 to Rutgers University for its support
of this research through its funding of the first author's visits to
Rutgers.}
 \and S. Shelah \\Hebrew University
and \\Rutgers University\thanks{ Partially supported by Basic Research Fund,
Israeli Academy of Sciences. Pub. No. 461}}

\begin{document}
\maketitle

\begin{abstract}
A new construction is given of non-standard uniserial modules over
certain valuation domains; the construction resembles that of a special
Aronszajn tree in set theory. A consequence is the proof of a sufficient
condition for the existence of non-standard uniserial modules; this
is a theorem
of ZFC which complements an earlier independence result.
\end{abstract}

\section*{Introduction}

This paper is a sequel to \cite{ES}. Both papers deal with the existence
of non-standard uniserial modules over valuation domains; we refer to
\cite{ES} for history and motivation. While the
main result of the previous paper was an independence result, the main
results of this one are theorems of ZFC, which complement and extend
the results of \cite{ES}.

We are interested in necessary and sufficient conditions for a valuation
domain $R$ to have the property that there is a non-standard uniserial
$R$-module of a given type $J/R$. (Precise definitions are given below.)
The question is interesting only when $R$ is uncountable, and since
additional complications arise for higher cardinals, we
confine ourselves to rings of \cha. Associated to any type $J/R$ is
an invariant, denoted $\gG(J/R)$, which is a member of a Boolean algebra
$D(\goo)$
(equal to ${\cal P}(\omega_1)$ modulo the filter of closed unbounded
sets). For example, if $R$ is an almost maximal valuation domain, then
$\gG(J/R) = 0$ for all types $J/R$; but there are natural and easily
defined examples where $\gG(J/R) = 1$.

 It is a fact that

\begin{quote}
{\em if $\gG(J/R) = 0$, then
there is no non-standard uniserial $R$-module of type $J/R$}
\end{quote}
 (cf.
\cite[Lemma 5]{ES}). In
\cite{ES} we showed that the converse is independent of ZFC + GCH;
the
consistency proof  that the converse fails involved the
construction of a valuation domain $R$ associated with a stationary and
co-stationary subset of $\goo$ --- that is, $0 < \gG(J/R) < 1$. The
existence of such sets 
requires a use of the Axiom of Choice; no such set can be explicitly
given. Thus --- without attempting to give a mathematical
definition of ``natural''  ---  we could say that for {\em natural}
 valuation
domains,
$R$, it is the case that for every type $J/R$,
$\gG(J/R)$ is either 0 or 1. For natural valuation domains,
it turns out that the converse {\it is}
true: if there is no non-standard uniserial $R$-module of type $J/R$,
then $\gG(J/R) =0$. This is a consequence of the following result which
is proved below (for {\em all} valuation domains of cardinality $\ha_1$):

\begin{quote}
{\it if $\gG(J/R) = 1$, then there is a non-standard uniserial
$R$-module of type $J/R$.}
\end{quote}

\noindent
(Theorem~\ref{G=1}.) This vindicates a conjecture made by Barbara Osofsky
 in \cite[(9), p. 164]{O1}. (See also the Remark following Theorem ~\ref{G=1}.)

The proof of Theorem~\ref{G=1} divides into several cases; the key new
 result which is used
 is a
construction of a non-standard uniserial module in the {\em essentially
countable} case; this construction is done in  ZFC
 and is motivated by
the construction of a special Aronszajn tree. (See Theorem~\ref{main}.)
Moreover, the
uniserial constructed is
``explicitly non-standard'' in that there is an associated
``special function'' which demonstrates that it is non-standard.
This special function continues to serve the same purpose in any
extension of the universe, V, of set theory,
so the module is ``absolutely'' non-standard.
In contrast, this may not be
the case with non-standard uniserials constructed using a prediction
(diamond) principle. (See the last section.)

The first author would like to thank L. Salce and S. Bazzoni for their
critical reading of a draft of this paper, and A. Kechris for a helpful
conversation on absoluteness.

\section*{Preliminaries}

For any ring $R$, we will use $R^*$ to denote the group of units of
$R$. If $r \in R$ we will write $x \equiv y \pmod{r}$ to mean $x - y
\in rR$.

A module is called uniserial if its submodules are linearly
ordered by inclusion. An integral domain $R$ is called a
valuation domain if it is a uniserial $R$-module.
 If $R$ is a valuation domain, let $Q$ denote its quotient
field; we assume $Q \neq R$. The residue field of $R$ is $R/P$, where
$P$ is the maximal ideal of $R$. \cite{FS} is a general reference for
modules over valuation domains.

If $J$ and $A$ are $R$-submodules of $Q$ with $A \subseteq J$, then
$J/A$ is a uniserial $R$-module, which is said to be  {\it standard}.
 A uniserial
$R$-module $U$ is said to be {\it non-standard} if it is not isomorphic
to a standard uniserial.

Given a uniserial module $U$, and a non-zero
element, $a$, of $U$, let $\ann(a) = \set{r \in R}{ra = 0}$ and
let $D(a) = \cup \{r^{-1}R \colon r$ divides a
in $U\}$. We say $U$ is  of {\em type} $J/A$ if $J/A \cong D(a)/\ann(a)$.
This is well-defined in that if $b$ is another non-zero element of $U$,
then $D(a)/\ann(a) \cong D(b)/\ann(b)$. For example,
$U$ has type $Q/R$
if and only if $U$ is divisible torsion and  the annihilator
ideal of every non-zero element of
$U$ is principal. (But notice that there is no $a \in U$ with $\ann(a) =
R$.) It is not hard to
see that if $U$ has type $J/A$, then
 $U$ is standard if and only if it is isomorphic to
$J/A$. We will only consider types of the form $J/R$; it is a
consequence of results of \cite{BFS} that the question of the
existence of a non-standard uniserial $R$-module of type $J_1/A$
can always be reduced to the question of  the existence of a
non-standard uniserial of type $J/R$
for an appropriate $J$.

>From now on we will assume that $R$ has cardinality $\hao$. We always have
$$
J = \cup _{\gs <\omega _1}r^{-1}_\gs R
\leqno(*)
$$
for some sequence of elements $\{r_\gs \colon \gs  < \omega _1\}$ such that for all
$\gt  < \gs $, $r_\gt |r_\gs $.  If $J$ is countably generated, then $U$ is
standard, so generally we will be assuming that $J$ is not countably generated;
then it has a set of generators as in ($\ast$), where, furthermore,
 $r_\gs $ does
not divide $r_\gt $ if $\gt  < \gs $.

If $\gd \in \lim(\goo)$, let

$$
J_\delta 
\deq \cup _{\gs <\delta } r^{-1}_\gs R.
\leqno(\ast \ast)
$$

By results in
[BS] every uniserial module $U$,
of type $J/R$, is described up to isomorphism by a family of units,
$\{e^\tau _\sigma \colon \sigma  < \tau  < \omega _1\}$ such that 
$$
e^\delta _\tau e^\tau _\sigma  \equiv e^\delta _\sigma  \pmod{r_\sigma}
\leqno(\dag)
$$
for all $\sigma  < \tau  < \delta  < \omega _1.$ 
Indeed, $U$ is a direct limit of  submodules $a_\gs R$
where $\ann(a_\gs) = r_\gs R$; then $a_\gs R
 \cong
r_\gs^{-1}R/R$ and $U$ is isomorphic to a direct
limit of the $r_\gs^{-1}R/R$, where the morphism from $r_\gs^{-1}R/R$
to $r_\gt^{-1}R/R$ takes $r_\gs^{-1}$ to $e_\gs^{\gt} r_\gs^{-1}$
if $a_\gs = e_\gs^\gt r_\gs^{-1} r_\gt a_\gt$.

If $U$ is given by (\dag ), then $U$ is standard if and only if there
exists a family $\{c_\sigma \colon \sigma  < \omega _1\}$ of units of $R$ such that
$$
c_\tau  \equiv e^\tau _\sigma c_\sigma  \pmod{r_\sigma}
\leqno{(\dag\dag)}
$$
for all $\sigma  < \tau  < \omega _1.$
Indeed,  if the
family $\{c_\sigma \colon \sigma  < \omega _1\}$ satisfying ($\dag\dag$)
exists, then multiplication by the $c_\gs$ gives rise to isomorphisms
from $r_\gs^{-1}R/R$ to $a_\gs R$, which induce an isomorphism of $J/R$
with $U$.




\section*{Essentially Countable Types}

\noindent
{\bf Definition.} Suppose $J = \cup _{\gs <\omega _1}r^{-1}_\gs R$ as
in (*).
Call the type $J/R$  {\it essentially uncountable} if
 for every $\gs  < \omega _1$ there exists $\gt  > \gs $
such that $r_\gs R/r_\gt R$ is uncountable. Otherwise, $J/R$ is {\it
essentially countable}; this is equivalent to saying that there is a
$\gamma < \go_1$ such that for all $ \gamma < \gs < \go_1$, $r_\gamma R/r_\gs R$
is countable.
Say that $J/R$ is {\it strongly countable} if for all $\gs < \go_1$,
$R/r_\gs R$ is countable; clearly, a strongly countable type is
essentially countable.

It is easily seen that the notions of being essentially or strongly
countable are well-defined, that is, independent of the choice of the
representation (*). If the residue field of $R$ is uncountable, then,
except in trivial cases, the types $J/R$ have to be essentially
uncountable; but if the residue field is countable, the question is more
delicate.

\begin{proposition}
\label{essun}
If the residue field of $R$ is uncountable, then every type $J/R$ 
such that $J$ is not countably generated is  essentially 
uncountable. 

\end{proposition} 

\proof 
Let $J = \cup _{\gs <\omega _1}r^{-1}_\gs R$ as in (*).
It suffices to prove that if $\gs  < \gt $, then $r_\gs R/r_\gt R$ is
uncountable. But $r_\gs R/r_\gt R \cong  R/tR$ where $t = r_\gt r^{-1}_\gs
 \in  P$. So we have $tR \subseteq  P \subseteq  R$, and hence
$(R/tR)/(P/tR) \cong  R/P$, the residue field of R. Since $R/P$ is uncountable, 
so is $R/tR$.  \fin 

\begin{theorem}
\label{exist} 

For any countable field $K$ there are valuation domains $R_1$ and
$R_2$, both  of cardinality $\hao$  with the same
residue field $K$ and   the same value group, whose quotient fields, $Q_1$
and $Q_2$, respectively, are generated by $\aleph _1$ but not countably many 
elements, and such that $Q_1/R_1$ is essentially uncountable and $Q_2/R_2$ is 
strongly countable.
\end{theorem}

\proof  
 Let $G$ be the ordered abelian group which is the direct
sum $\oplus _{\alpha <\omega _1}{\open Z}\alpha $ ordered 
anti-lexicographically; that is, $\Sigma _\alpha  n_\alpha \alpha  > 0$ if and 
only if $n_\beta  > 0$, where $\beta $ is maximal such that $n_\beta  \neq  0$.
 In particular, the 
basis elements have their natural order and if $\alpha  < \beta $, then 
$k\alpha  < \beta $ in $G$ for any $k \in  {\open Z}.$  Let $G^+ = \{g \in  G\colon
g \geq  0\}$.

Let
$\hat{R}  = K[[G]]$, that is, $\hat{R} = \{\sum_{g\in \Delta } k_gX^g\colon k_g \in  K, \Delta \hbox{ a 
well-ordered subset of } G^+\}$, with the obvious addition and
multiplication (cf. \cite[p. 156]{O1}).
Given an element $y = \sum _{g\in \Delta } k_gX^g$ of $\hat{R} $, let 
$\supp(y) = \{g \in  \Delta \colon k_g \neq  0\}$; let $\psupp(y) = \{\alpha  \in  
\omega _1\colon \exists g \in  \supp(y)$ whose projection on ${\open Z}\alpha $ is 
non-zero$\}$. Define $v(y) = $ the least element of $\supp(y)$. If $X \subseteq  G$, then $y\rest X$ is defined to be 
$\sum _{g\in X\cap \Delta } k_gX^g$. Let $y\rest \nu  = y\rest \{g \in  G^+\colon g 
< \nu \}.$ 

Let $R_1 = \{y \in  \hat{R} \colon p$-supp$(y)$ is finite$\}$. Then $R_1$ is a valuation
domain since $p$-supp$(xy^{-1}) \subseteq  p$-supp$(x) \cup  p$-supp$(y)$. Let 
$R_2$ be the valuation subring of $R_1$ generated by $\{X^g: g \in  G\}$. We
have $Q_j = \cup _{\alpha  < \omega _1} X^{-\alpha }R_j$ for $j = 1$, 2. Now 
$Q_1/R_1$ is essentially uncountable since for all $\beta  > \alpha $, 
$X^\alpha R_1/X^\beta R_1$ contains the $2^{\aleph _0}$ elements of the form
$$
\sum _{n \in  \omega } \zeta (n)X^{(n + 1)\alpha }
$$

\noindent
(with $p$-supp $= \{\alpha \})$ where $\zeta $ is any function: $\omega  
\rightarrow  2$. $R_1$ has cardinality $2^{\aleph _0}$; if $2^{\aleph _0} > 
\aleph _1$, to get an example of cardinality $\aleph _1$, choose a 
valuation subring of $R_1$ which contains all the monomials $X^g$ ($g \in  G)$ 
and $\aleph _1$ of the elements $\sum _{n \in  \omega } 
\zeta (n)X^{(n + 1)\alpha }$ for each $\alpha .$ 

We claim that $Q_2/R_2$ is essentially countable. Let $K[G]$ be the subring of 
$\hat{R}$ generated by $\{X^g \colon g \in  G^+\}$; thus $K[G]$ consists of the elements of
$\hat{R}$ with finite support; we shall refer to them as {\it polynomials}. $R_2$ 
consists of all elements of the form $xy^{-1}$ where $x$ and $y$ are 
polynomials and $v(x) \geq  v(y)$.
 We claim
that $ R_2/X^\beta R_2$ is countable for any $\beta < \goo$.
 There are uncountably many polynomials, but
we have to show that there are only countably many truncations 
$xy^{-1}\rest \beta $.

Given polynomials $x$ and $y$ with $v(x) \geq  v(y)$, there is a finitely 
generated subgroup $\oplus _{1\leq i\leq d} {\open Z}\sigma _i$ of $G$ (with
$\sigma _1 < \sigma _2 <  \ldots < \sigma _d)$ such that $x$ and $y$ are linear
combinations of monomials $X^g$ with $g \in  \oplus _{1\leq i\leq d}
{\open Z}\sigma _i$. More precisely, there exist $k, r \in \go$ and a $(k + r)$-tuple $(a_1$,
$ \ldots ,a_{k+r})$ of elements of $K$ and $k + r$ linear terms $t_j$ of the form
$$
t_j = \sum_{i = 1}^d n_{ij}v_i
$$

\noindent
$(n_{ij} \in  {\open Z}$, $v_i$ variables) such that if we let $t_j(\sigma )$ 
denote $\sum_{i = 1}^d n_{ij}\gs_i$, then

\labpar{\natural}{
$x = \sum_{j = 1}^k a_jX^{t_j(\sigma )}\hbox{ and }\ \ y = 
\sum_{j = k+1}^r a_jX^{t_j(\sigma )}.$}

\noindent
Finally, there is $q \leq  d$ such that $\sigma _q$ is maximal with $\sigma _q
< \beta.$

Now, consideration of the algorithm for computing $xy^{-1}$ shows that,
for fixed $(a_1, \ldots, a_{k + r})$ and $t_j$,  there 
are linear terms 
$$
s_\ell  = \sum_{i = 1}^d m_{i\ell }v_i
$$

\noindent
$(m_{i\ell } \in  {\open Z}$, $\ell  \in  \omega )$ and elements $c_\ell  \in  
K$ such that for {\it any} strictly increasing sequence $\sigma  = 
\langle \sigma _1,  \ldots , \sigma _d\rangle $, if $x$ and $y$ are as in ($\natural$),
then
$$
xy^{-1} = \sum _{\ell  \in  \omega } c_\ell X^{s_\ell (\sigma )}.
$$

\noindent
For any $q \leq  d$, only certain of the $s_\ell $ involve only variables $v_i$
with $i \leq  q$ (i.e. $m_{i\ell } = 0$ if $i > q$); say these are the
$s_\ell$ with $\ell \in T$ ($T \subseteq \go$). If $\gs$ is such that
$\gs_i < \gb $ iff $i \leq q$, then
 $ xy^{-1}\rest\gb  = \sum _{\ell  \in  T } c_\ell X^{s_\ell
(\sigma )}$.

There are only countably many choices for $q$, $d$, $k$ and $r$ in $\go$, $(a_1$,  \ldots ,
$a_{k+r}) \in  K^{k+r}$, and for $\sigma _1 < ..$. $< \sigma _q < \beta
$. Therefore, there are only countably many possibilities for the
truncations $xy^{-1}\rest \beta $.  \fin

By the first part of the following,
 the type $Q_2/R_2$ of the previous theorem {\em must} be strongly
countable; on the other hand, there are types which are essentially
countable but not strongly countable.

\begin {proposition}
\label{strgly} 

(i) If $Q/R$ is essentially countable, then it is strongly countable. 

(ii) For any countable field $K$, there is a valuation domain $R$ with residue 
field $K$ which has a type $J/R$ which is essentially countable but not 
strongly countable. 

\end {proposition} 

\proof 
(i) Since $Q/R$ is essentially countable, we can write $Q = 
\cup _{\gs <\omega _1} r^{-1}_\gs R$ such that for all $\gs  < \gt $,
$r_\gs R/r_\gt R$ is countable. We claim that $R/r_0R$ is countable, which
clearly is equivalent to $R/r_\gt R$  countable for all $\gt  < \omega _1$.
Suppose not. There is a $\gs  < \omega _1$ such that $r_\gs  = r^2_0t$ for some
$t \in  R$ (since $r^{-2}_0 \in  Q)$. But then $r_0R/r_\gs R \cong  R/tr_0R$,
which is uncountable since $R/r_0R$ is uncountable, and this contradicts the 
choice of the $r_\gs .$

(ii) Let $G = \oplus _{\alpha \leq \omega _1} {\open Z}\alpha $, ordered 
anti-lexicographically. Let $\hat{R} =K[[G]]$ (cf. proof of Theorem~\ref{exist}), and let $R$ be
the smallest valuation subring of $\hat{R}$ containing all the monomials $X^g\ ( g 
\in  G)$. Let $J = \cup _{\alpha <\omega _1} r^{-1}_\alpha R$ where $r_\alpha  
= X^{\alpha +\omega _1}$. Then the proof that $r_\alpha R/r_\beta R \cong  
R/X^{\beta -\alpha }R$ is countable for all $\ga < \gb$ is the same as in Theorem~\ref{exist}. But $R/r_0R =
R/X^{\omega _1}R$ is clearly uncountable.  \fin

\noindent
{\bf Remark. } \ More generally, referring to a dichotomy in \cite[Prop.
7, p. 155]{O1}, if the type $J/R$ is essentially countable and falls
into Case (A), then $J/R$ is strongly countable; if it falls into case
(B), then it is not strongly countable.

\section*{Gamma Invariants}

 A subset $C$ of $\go_1$ is called a {\em cub} ---
short for closed unbounded set --- if $\sup C = \go_1$ and for all
 $Y \subseteq  C$,
$\sup Y \in  \go_1 $ implies  $\sup Y \in  C$.
 Call two subsets, $S_1$ and $S_2$, of $\go_1$
equivalent iff there is a cub $C$ such that $S_1 \cap C = S_2 \cap C$.
Let $\tilde{S}$ denote the equivalence class of $S$. The inclusion
relation induces a partial
order on the set, $D(\goo)$, of equivalence classes, i.e.,
$\tilde{S_1} \leq \tilde{S_2}$ if and only if there is a cub $C$ such
that $S_1 \cap C \se S_2 \cap C$. In fact, this induces a Boolean
algebra structure on $D(\goo)$, with least element, $0$, the equivalence
class of sets disjoint from a cub; and greatest element, $1$, the
equivalence class of sets containing a cub. We say $S$ is {\em
stationary} if $\tilde{S} \neq 0$, i.e., for every cub $C$, $C \cap S
\neq \emptyset$. We say $S$ is {\em co-stationary} if $\goo \setminus S$
is stationary.

Given $R$ and a type $J/R$, where $J$ is as in (*), define $\Gamma(J/R)$ to be
$\tilde{S}$, where

$$  S = \{\delta  \in  \lim (\omega _1)\colon
 R/\cap _{\gs <\delta }r_\gs R\hbox{ is not complete}\}
$$

\noindent
where the topology on $R/\cap _{\gs <\delta }r_\gs R$ is the metrizable linear topology
with a basis of neighborhoods of 0 given by the submodules $r_\gs R\  (\gs  <
\delta ).$ This definition is independent of the choice of the
representation of $J$ as in (*) --- see \cite{ES}.

For any limit ordinal $\delta  < \omega _1$, let

$${\cal T}^\delta _{J/R} = \{\langle u_\sigma \colon \sigma  < \delta \rangle \colon
\forall \sigma  < \tau  < \delta (u_\sigma  \in  R^{*},\  \hbox{and } u_\tau  - u_\sigma 
\in  r_\sigma R)\};$$

\noindent
that is, ${\cal T}^\delta _{J/R}$ consists of sequences of units which are
Cauchy in the metrizable topology on $R/\cap _{\gs <\delta }r_\gs R$. Let ${\cal L}^\delta _{J/R}$ consist of those
members of ${\cal T}^\delta _{J/R}$ which have limits in $R$, i.e.

$${\cal L}^\delta _{J/R}= \{\langle u_\sigma \colon \sigma  < \delta \rangle  \in
{\cal T}^\delta _{J/R} \colon \exists u_\delta  \in  R^{*} \hbox{
s.t. }\forall \sigma  < 
\delta (u_\delta  - u_\sigma  \in  r_\sigma R)\}.$$
Note that $\Gamma(J/R) = \tilde{S}$ where
$$
S = \{\delta  \in  \lim (\omega _1)\colon {\cal T}^\delta _{J/R} \neq
{\cal L}^\delta _{J/R}\}.
$$

If $J$ is not countably generated, then $\neg$CH implies
that $\gG(J/R) = 1$, since the completion of   $ R/\cap _{\gs <\delta
}r_\gs R$ has cardinality $2^{\ha_0} > \hao$.
An $\omega _1${\it -filtration of $R$ by subrings} is
an increasing chain $\{N_\alpha \colon \alpha  \in  \omega _1\}$ of
countable subrings of $R$ such that $R = \cup _{\alpha \in \omega _1}
N_\alpha $, and for limit $\alpha $, $N_\alpha  = \cup _{\beta <\alpha }
N_\beta $.

Define  $\Gamma'(J/R) =
\tilde{E}'$ where  
\begin{quotation}
\noindent
$E' =
\{\delta  \in  \lim (\omega _1)\colon \exists \langle u_\sigma \colon \sigma  < 
\delta \rangle  \in  {\cal T}^\delta _{J/R} \hbox{ s.t. } \forall f \in  R^{*} \exists
\sigma  < \delta \ \hbox{ s.t. } \\
u_\sigma f \notin  N_\delta \pmod{r_\sigma} \}.$
\end{quotation}

Again, it can be shown that the definition does not
depend on the choice of $\{r_\nu \colon \nu  < \omega _1\}$ or of $\{N_\alpha \colon 
\alpha  < \omega _1\}$.
 Notice that $\Gamma'(J/R) \leq  \Gamma(J/R)$
since if ${\cal T}^\delta _{J/R} = {\cal L}^\delta _{J/R}$, then we can let
$f$ be a limit of $\langle u_\gs^{-1} \colon \gs < \gd \rangle$.

In \cite[Theorem 7]{ES} it is proved that if $\Gamma'(J/R) \neq  0$, then there is a
non-standard uniserial $R$-module.

\begin{theorem} 
\label{gamma}
Suppose $J/R$ is essentially countable. Then 

(i) $\Gamma '(J/R) = 0;$

(ii) $\Gamma (J/R) = 1.$ 
\medskip
\end{theorem} 
\medskip
\proof
Without loss of generality we can assume that $J = \cup _{\gs  < \omega _1}
r^{-1}_\gs R$ where $r_0R/r_\gs R$ is countable for all $\gs  < \omega _1$.

(i) We can also assume that the $\omega _1$-filtration of $R$ by subrings, $R = 
\cup _{\alpha  < \omega _1} N_\alpha $, has the property that for all 
$\alpha $, $N_\alpha $ contains a complete set of representatives of 
$r_0R/r_\sigma R$ for each $\sigma  < \ga $. For any $\delta  \in  
\lim (\omega _1)$,  and any $\langle u_\sigma : \sigma  < \delta \rangle $ in 
${\cal T}^\delta _{J/R}$, let $f = u^{-1}_0$. To show that $\delta
\notin  E'$,
 it suffices to show that $u_\sigma f
\in  N_\delta $ (mod $r_\sigma)$ for all $\sigma  < \delta $. Now $u_\sigma f
= u_\sigma u^{-1}_0 \equiv  1$ (mod $r_0)$, since $u_\sigma  \equiv  u_0$ (mod
$r_0)$, by definition of ${\cal T}^\delta _{J/R}$. Say $u_\sigma u^{-1}_0 - 1 
= y \in  r_0R$. By the assumption on $N_\delta $, there exists $a \in  
N_\delta $ such that $y \equiv  a$  (mod $r_\sigma)$. Then $u_\sigma f = 
u_\sigma u^{-1}_0 = 1 + y \equiv  1 + a$ (mod $r_\sigma)$, and $1 + a \in  
N_\delta $ since $N_\delta $ is a subring of $R.$ 

(ii) To show that $\Gamma (J/R) = 1$, it suffices to show that for all limit 
ordinals $\delta  < \omega _1$, $R/\cap _{\nu  < \delta } r_\nu R$ is not
complete. Assuming that it is complete, we shall obtain a contradiction by 
showing that $r_0R/r_\delta R$ is uncountable. Fix a ladder on $\delta $, 
i.e., a strictly increasing sequence $\langle \nu _n: n \in  \omega \rangle $
whose sup is $\delta $. For each function $\zeta \colon \omega  \rightarrow  2$,
define $u^\zeta  = \langle u^\zeta _\sigma : \sigma  < \delta \rangle  \in  
{\cal T}^\delta _{J/R}$ as follows: if $\nu _m < \sigma  \leq  \nu _{m+1}$,
then 
$$
u_\sigma  = \sum_{i\leq m} \zeta (i)r_{\nu _i}.
$$

\noindent
Clearly $u_\sigma  \in  r_0R$, and if $\tau  > \sigma $, where $\nu _k < \tau  
\leq  \nu _{k+1}$, then $m \leq  k$ and 
$$
u_\tau  - u_\sigma  = \sum_{i = m+1}^k \zeta (i)r_{\nu _i} \in  
r_{\nu _{m+1}}R \subseteq  r_\sigma R.
$$

\noindent
Since $R/\cap _{\nu  < \delta } r_\nu R$ is assumed to be complete,  for each 
$\zeta $ there is an element $u^\zeta _* \in  R$ which represents the limit of 
$\langle u^\zeta _\sigma : \sigma  < \delta \rangle $ in $R/\cap _{\nu  < 
\delta } r_\nu R$. To obtain a contradiction, we need only show that if $\eta  
\neq  \zeta $, then $u^\zeta _* - u^\eta _* 
\notin  r_\delta R$. Without loss of generality there exists $m$  such that $\zeta \rest m = \eta \rest m$ and
$\zeta (m) = 0$, $\eta (m) = 1$, then 
$$
u^\eta _* - u^\eta _{\nu _{m+1}} \in  r_{\nu _{m+1}}R\hbox{;  and}
$$
$$
u^\zeta _* - u^\zeta _{\nu _{m+1}} \in  r_{\nu _{m+1}}R;
$$

\noindent
but $u^\eta _{\nu _{m+1}} - u^\zeta _{\nu _{m+1}} = r_{\nu _m} \notin  
r_{\nu _{m+1}}R$, so $u^\zeta _* - u^\eta _* \notin  r_{\nu _{m+1}}R \supseteq 
r_\delta R$.  \fin 



\section*{Special Aronszajn trees}

This section contains standard material from set theory.
(See, for example, \cite[\S 22]{J} or \cite[Ch.
7, \S 3]{Dr}.) It is included simply to provide
motivation for the notation and proof  in the next section.

A {\it tree} is a partially ordered set $(T$, $<)$
such that the
predecessors of any element are well ordered. An element $x$ of $T$ is
 said to have {\it height} $\alpha $, denoted  
${\rm ht}(x) = \alpha $, if the order-type of $\{y \in  T\colon y < x\}$
 is $\alpha $. The
{\it height} of $T$ is defined to be
$\sup \{{\rm ht}(x) + 1\colon x \in  T\}.$ 
If $T$ is a tree, a {\it branch}  of $T$ is a maximal linearly ordered 
initial subset of $T$; the {\it length} of a branch is its order type.
  If $T$ is a tree, let $T_\alpha  = \{y \in  T: ht(y) = \alpha \}$.
We say that a tree $T$ is a {\em $\gk$-Aronszajn tree} if $T$ is of
height $\gk$,  $|T_\ga| < \gk$ for every $\ga < \gk$, and $T$ has no
branch of length $\gk$. 

 A tree $T$ of 
height $\omega _1$ is a {\it special Aronszajn tree} if $T_\alpha $ is 
countable for all $\alpha  < \omega _1$ and for each $\alpha  < \omega _1$ 
there is a function $f_\alpha \colon T_\alpha  \rightarrow  {\open Q}$ such that
\labpar{\$}{
whenever $x \in  T_\alpha $ and $y \in  T_\beta $ and $x < y$, then \\
$f_\alpha (x) < f_\beta (y)$.}
 Notice that a special Aronszajn tree
is an $\go_1$-Aronszajn tree,
 since an uncountable branch would give rise to an uncountable 
increasing sequence of rationals. 

K\"onig's Lemma implies that
 there is no $\go$-Aronszajn tree.
However, there is an $\go_1$-Aronszajn tree:

\begin{theorem} 
\label{artree}
There is a special Aronszajn tree. 
\end{theorem} 

\proof 
Let $^{<\alpha }\omega $ denote the set of all functions from $\set{\gb
\in \goo}{\gb < \ga}$ to $\go$.
We shall construct $T_\alpha $ and $f_\alpha $ by induction on $\alpha  <
\omega _1$ such that $T_\alpha $ is a countable subset of $^{<\alpha }\omega $
and the partial ordering is inclusion, i.e., if $x \in  T_\alpha $ and $y 
\in  T_\beta $ then $x < y$ if and only if $\alpha  < \beta $ and 
$y\rest \alpha  = x.$ Finally, $T$ will be defined to be
$\bigcup_{\ga < \go_1} T_\ga$.

Let $T_0 = \{\emptyset \}$, $f_0(\emptyset) = 0$,  $T_1 = {}^{\{0\}}\omega $, and $f_1: T_1 \rightarrow
{\open Q}$ be onto $(0, \infty)$. Suppose now that $T_\alpha $ and $f_\alpha $
have been defined for all $\alpha  < \gd $ such that for all $\gs < \gr
< \gd$:

 \labpar{\star}{for any $\epsilon > 0$, and $x \in T_\gs$
 there is $y \in  T_\gr $ such that $x < y$ and
  $f_\gr (y) <
f_\gs (x) + \epsilon .$}

There are two cases. In the first case, if $\gd $ is a successor ordinal, 
$\gd  = \gt  + 1$, let
$$
T_\gd  = \{x \cup  \{(\gt \hbox{, } n)\} : n \in  \omega \hbox{, } x \in 
T_\gt  \}.
$$

\noindent
Define $f_\gd $ so that for every $x \in  T_\gt ,$
$$
\{f_\gd (x \cup  \{(\gt \hbox{, } n)\}): n \in  \omega \} = \{r \in  
{\open Q}: r > f_\gt (x)\}.
$$

\noindent
Clearly ($\star$) continues to hold.  

In the second case, $\gd $ is a limit ordinal. Choose a ladder
$\langle \gn_n: n \in  \omega \rangle $ on $\gd $. For each
$\gs  < \gd $, $x \in  T_\gs $ and $k > 0$, by inductive hypothesis ($\star$)
there exists a sequence $\langle y_n: n \in  \omega$  s.t. $\gn_n >
\gs \rangle $ such that $y_n \in T_{\gn_n}$, $x < y_n < y_m$ for $n < m$ and $f_{\gn_n}(y_n) <
f_\gs (x) + (1/k - 1/n)$. Let $y[\gs $, $x$, $k] = \bigcup _{n\in \omega } y_n
\in  {}^{<\gd}\omega $ . Let $T_\gd $ consist of one such $y[\gs $, $x$,
$k]$ for each $\gs $, $x$, $k$. Define $f_\gd (y[\gs $, $x$, $k]) = 
f_\gs (x) + 1/k$. Then it is clear that ($\star$) still holds.  \fin

\section*{Special Uniserial Modules}

\noindent
{\bf Definition.}\ \ 
Suppose $U$ is a uniserial module of type $J/R$ where $J = \cup _{\gs  < 
\omega _1} r^{-1}_\gs R$ as in (*). For each $\gs > \go_1$, fix an element $a_\gs $ of $U$ such that 
$\ann(a_\gs ) = r_\gs R$ (so that the submodule $a_\gs R  $ of $U$ is
isomorphic to  $R/r_\gs R \cong  
r^{-1}_\gs R/R)$. Let $I_\gs $ be the set of all $R$-module isomorphisms
$\varphi \colon
a_\gs R \rightarrow r^{-1}_\gs R/R$. We say that $\set{f_\gs}{\gs \in \goo}$
is a {\em special family of functions} for  $U$
 if for each $\gs  < \omega _1$,  $f_\gs :
I_\gs  \rightarrow  {\open Q}$ such that whenever $\gs  < \gr $ and $\varphi  
\in  I_\gr $ extends $\psi  \in  I_\gs $, then $f_\gs (\psi ) < 
f_\gr (\varphi )$.

\begin{lemma}
\label{special}

If $U$ has a special family of functions, then $U$ is non-standard.

\end{lemma} 

\proof
Suppose there is an isomorphism $\theta \colon U \rightarrow  J/R$. Then for every
$\gs  < \omega _1$, $\theta $ restricts to an isomorphism $\varphi_\gs$
of $a_\gs R $ onto 
$r^{-1}_\gs R/R$. But then $\langle f_\gs(\varphi_\gs) \colon \gs < \goo
\rangle$ is an uncountable
strictly increasing sequence of rationals, a contradiction.  \fin

With this lemma as justification, we will     say that $U$ is {\it
explicitly non-standard} if $U$ has a special family of functions.

If the uniserial module uniserial module $U$, 
of type $J/R$, is described up to isomorphism by a family of units, 
$\{e^\gr _\sigma \colon \sigma  < \gr  < \omega _1\}$ as in 
(\dag), 
then it is clear that $U$ is explicitly non-standard if and only if for every $\gs <
\go_1$, there is a function $$f_\gs \colon (R/r_\gs R)^\ast
\longrightarrow \oQ$$
such that

\labpar{\$\$}{whenever $\gs < \gr$ and $c_\gs, c_\gr \in R^\ast$ satisfy $c_\gr
\equiv c_\gs e_\gs^\gr  \pmod{r_\gs}$, then $f_\gs(c_\gs) < f_\gr(c_\gr)$.}

\noindent
(Here, and hereafter,  we abuse notation and regard $f_\gr$ and $f_\gs$ as functions
on $R^\ast$.)

           Note that we have a tree, T, such that $T_\gs =
(R/r_\gs R)^\ast$ and the partial ordering is given by: $$c_\gs + r_\gs R
< c_\gr + r_\gr R \  \Longleftrightarrow \
\gs < \gr \hbox{ and  } c_\gr \equiv c_\gs e_\gs^\gr
\pmod{r_\gs}.$$
Assume $\gs < \gr$.
Each $c_\gs $ has at least one successor of height $\gr$, namely $c_\gs
e_\gs^\gr$, and if
 $r_\gs R/r_\gr R$ is countable, then  $c_\gs$ has
only countably many successors of height $\gr$.
For each $c_\gr \in  T_\gr$, its unique
predecessor in $T_\gs$ is $c_\gr (e_\gs^\gr)^{-1} $.
(Here again we abuse notation and write, for example, $c_\gr$ for an
element of $T_\gr$ instead of $c_\gr + r_\gr R$.)

Without loss of generality we can assume that the $r_\gs$ are such that for all
$\gs < \goo$, $r_\gs R/r_{\gs + 1} R$ is infinite. (Just choose a
subsequence of the original $r_\gs$'s if necessary.) Thus for all $\gs <
\goo$,
 there is an infinite subset $W_\gs$ of $R^*$
such that for all $u \neq v \in W_\gs$, $u \equiv 1 \pmod{r_\gs}$ and
$u \not\equiv v \pmod{r_{\gs + 1}}$.

\begin{theorem} 
\label{main} 

If $J/R$ is an essentially countable type, then there is an explicitly 
non-standard uniserial $R$-module of type $J/R.$ 

\end{theorem} 

\proof 
We will first give the construction in the case when $J/R$ is strongly 
countable, and afterward indicate the modifications needed for the general 
case. Thus $T_\gs$ is assumed countable for all $\gs < \goo$.

We will define, by induction on $\delta $, $e^\tau _\sigma $ for $\sigma  <
\tau  < \delta $ as in (\dag ) and, at the same time, the maps $f_\sigma \colon
(R/r_\sigma R)^\ast \rightarrow  {\open Q}$ for $\sigma  < \delta $. We will do this
so that (\$\$) holds and the following condition is satisfied for all
$\gs < \gr < \goo$:

\labpar{\star\star_{\gs, \gr}}{for any $\epsilon  > 0$,  $m \in  \omega $,
$c_\gs^j \in T_\gs$,
and $c^j_\rho  \in  T_\gr$    ($j = 1, \ldots , m)$ such that
$c_\gs^j < c_\gr^j$, there exists $u
\in  R^\ast$ such that $u \equiv  1 \pmod{r_\sigma}$ and $f_\rho (uc^j_\rho ) < 
f_\sigma (c^j_\gs) + \epsilon $ for all $ j = 1, \ldots, m$.}

Note that
the $c_\gr^j$ determine the $c_\gs^j$ --- $c_\gs^j \equiv
c_\gr^j(e_\gs^\gr)^{-1} \pmod{r_\gs}$ --- and $uc_\gr^j$ is another
successor of $c_\gs^j$ of height $\gr$.
For any given $\epsilon  > 0$, $\sigma  < \rho  < \delta $, $m \in  \omega $,
and $c^j_\rho  \in  R^\ast$    ($ j = 1, \ldots, m)$, there exist
infinitely many $u$ as in ($\star\star_{\gs, \gr}$), since we can decrease $\epsilon$
as much as we like.

Suppose we have
defined $e^\rho _\sigma $ and $f_\sigma $ for all $\sigma  < \rho  < \delta $
satisfying the inductive hypotheses.
 Let
$\langle u_n: n \in  \omega \rangle $ enumerate representatives of all the
elements of $(R/r_\delta R)^\ast $. Also, let
$\langle \gth_q: q \in  \omega \rangle $ enumerate all instances of
($\star\star_{\gs, \gd}$), for all $\gs < \gd$, with each instance repeated infinitely often.
More precisely, we enumerate (with infinite repetition) all
tuples of the form
$$
\langle \epsilon = \frac{1}{n},\  \gs, \   c^j_\gd + r_\gd R \colon j =
1, \ldots, m 
\rangle
$$
with $n \in \go \setminus \{0 \}$, $\gs < \gd$, and $c^j_\gd \in R^\ast$.

We will define $f_\gd$ as the union of a
 chain of functions $f_{\gd ,k}$ into
${\open Q}$, each with a finite domain. When $k$ is even we will
concentrate on insuring that the domain of $f_\gd$ will be $T_\gd$; and
when $k$ is odd, we will work at satisfying the conditions
($\star\star_{\gs, \gd}$).

Suppose first that $\delta  = \tau  + 1$ and
define $e^\delta _\tau  = 1$ and $e^\delta _\sigma  = e^\tau _\sigma $ for
$\sigma  < \tau $. Suppose that $f_{\gd, i}$ has been defined for $i <
k$, and assume first that $k$ is even. Let $n$ be minimal such that $u_n
\notin \dom(f_{\gd, k - 1})$. Let $\dom(f_{\gd, k}) = \dom(f_{\gd, k -
1}) \cup \{ u_n \}$ and  let $f_{\gd, k}(u_n)$ be any rational
greater than
$f_\gt(u_n(e_\gt^\gd)^{-1})$ ($ =f_\gt(u_n)$).

Now suppose $k$ is odd; say $k = 2q + 1$.
 It's
easy to see that it's enough to construct $f_\delta $ to satisfy
($\star\star_{\gt, \gd}$).
So if $\gth_q$ is an instance of ($\star\star_{\gs, \gd}$) for $\gs <
\gt$, let $f_{\gd, k} = f_{\gd, k - 1}$. Otherwise, suppose $\gth_q$ is the
instance of ($\star\star_{\gt, \gd}$) given by
$$ \frac{1}{n},\  c^j_\gd + r_\gd R \colon  j =
1, \ldots, m 
.$$
Since $W_\gt$ is infinite (see above),  there is a unit $u$ such that
$u \equiv  1 \pmod{r_\tau} $ and $uc^j_\delta  \notin  \dom(f_{\delta
,k-1})$
for $j = 1, \ldots, m$. Then define $f_{\delta ,k}$ to be the extension of
$f_{\delta ,k-1}$ with domain $= \dom(f_{\delta ,k-1}) \cup
\{uc^j_\delta  \colon  j 
= 1$, ..., $m\}$ such that $$f_{\delta ,k}(uc^j_\delta ) =
f_\gt(c^j_\gd (e_\gt^\gd)^{-1}) + \frac{1}{2n}.$$

Now we consider the case where $\delta $ is a limit ordinal.
Fix a ladder $\langle \nu _n: n \in  \omega \rangle $ on $\delta $.
We are going to define units $e_{\nu _n}^\gd$
  by induction  such that $e_{\nu _n}^\gd
\equiv  e_{\nu _m}^\gd e^{\nu _m}_{\nu _n} \pmod{r_{\nu _n}}$ whenever $n < m <
\omega $. This will easily determine the sequence $\langle e_\sigma^\gd
 : \sigma  <
\delta \rangle $ such that for all $\sigma  < \tau  < \delta $,
$e_\gs^\gd
\equiv  e_\gt^\gd e^\gt _\sigma  \pmod{r_\sigma} $;
then (\dag ) will be  satisfied for
$\langle e^\tau _\sigma : \sigma  < \tau  \leq  \delta \rangle $. 

\newcommand{\en}{e_n}
\newcommand{\ek}{e_k}
For simplicity of notation, let $\en$ denote $e_{\nu_n}^\gd$.
Suppose we've already defined $f_{\gd, k-1}$ and $\ek$ such that
for all $x \in \dom(f_{\gd, k-1})$,
$$ f_{\gn_k}(x\ek^{-1}) < f_{\gd, k-1}(x).$$
(Recall that if $x \in T_\gd$, then $x\en^{-1}$ is the unique predecessor
of x in $T_{\gn_n}$.) If $k$ is even, we proceed as in the even case
above (when $\gd$ is a successor). If $k = 2q + 1$ and $\gth_q$ is
$$
\langle  \frac{1}{n},\  \gs, \   c^j_\gd + r_\gd R \colon j =
1, \ldots, m 
\rangle
$$
we can assume --- since each instance is repeated infinitely often ---
that $\gs < \gn_k$. Thus $e_\gs^\gd = e_{\gn_k}^\gd e_\gs^{\gn_k}$ is
defined.
 Note that
$c^j_\delta (e_\gs^\gd)^{-1}  \equiv
c^j_\delta \ek^{-1}(e^{\nu_k}_\gs )^{-1} \pmod{r_\gs}$,
 so we can apply ($\star\star_{\gs, \gn_k}$)
[with  $c_{\gn_k}^j = c_\gd^j\ek^{-1}$] and
obtain a unit $w \equiv  1 \pmod{r_\gs} $ such that
for all $j = 1, \ldots, m$
$$
f_{\nu _{k}}(wc^j_\delta\ek^{-1}) < f_\gs (c^j_\delta (e_\gs^\gd)^{-1}) +
\frac{1}{2n}.
$$
Moreover, since there are infinitely many such $w$, we can choose one so that 
the elements $wc^j_\delta$ $(j = 1, \ldots, m)$
do not belong to $\dom(f_{\gd, k-1})$. Let these be the new
elements of the domain of $f_{\gd, k}$ and define
$$ f_{\gd, k}(wc_\gd^j) = f_\gs (c^j_\delta (e_\gs^\gd)^{-1}) +
\frac{1}{2n}.$$

Now we will define $e_{k+1}$ (for $k$ odd or even).
For each $x \in \dom(f_{\gd, k})$ we have committed ourselves to
$f_\gd(x)$ ($= f_{\gd, k}(x)$) and to the predecessor of $x$ in
$T_{\gn_k}$ ($= x\ek^{-1}$); we need to choose $e_{k+1}$ so that $x$
and its predecessor, $xe_{k+1}^{-1}$, in $T_{\gn_{k+1}}$ satisfy (\$\$).

 Let $e' =
e_k(e_{\nu _k}^{\nu _{k+1}})^{-1}$. The desired element $e_{k+1}$ will
have the form $ue'$ for some unit $u \equiv  1 \pmod{r_{\nu
_k}}$.
Choose $\epsilon' < f_{\gd, k}(x) - f_{\gn_k}(x\ek^{-1})$ for each $x \in
\dom(f_{\gd, k})$. Apply ($\star\star_{\gn_k, \gn_{k+1}}$) to this
$\epsilon'$ and   $x\ek^{-1} \in T_{\gn_k}$,
 $xe'^{-1} \in T_{\gn_{k+1}}$
 ($x \in \dom(f_{\gd, k})$).
(Note that $  x\ek^{-1} < xe'^{-1}$ by choice of $e'$.)
This gives us $v \equiv 1 \pmod{r_{\gn_k}}$ such that for all $x$
$$f_{\gn_{k+1}}(vxe'^{-1}) < f_{\gn_k}(x\ek^{-1}) + \epsilon' < f_{\gd,
k}(x).$$
Then we let $e_{k+1} = v^{-1}e'$, and we have completed the inductive
step.

\bigskip

This completes the proof in the strongly countable case. We turn now
to the general (essentially countable) case. In this case, $R/r_0R$
may be uncountable; let $Z$ be a complete set of representatives of $(R/r_0R)^\ast$.
Fix $z_0 \in Z$. We first define, by induction on $\gs$,
 $f_\gs(c_\gs)$ --- or, more precisely, $f_\gs(c_\gs +
 r_\gs R)$ ---  for all $c_\gs \in R^\ast$ such that $c_\gs \equiv
 z_0e_0^\gs
 \pmod{r_0}$. We do the construction  exactly as in
 the previous strongly countable case; this will work since there are
 only countably many cosets $c + r_\gs R$ such that $c \equiv z_0e_0^\gs
 \pmod{r_0}$ since $r_0R/r_\gs R$ is countable.

Having done this, the $e_\gs^\gt$
are determined.
We claim that there is no family $\{c_\gs \colon \gs < \goo \}$
satisfying ($\dag\dag$). Indeed, suppose we had such a family. Let $z
\in Z$ be such that $c_0 \equiv z \pmod{r_0}$. Then for all $\gs < \goo$,
$c_\gs \equiv ze_0^\gs \pmod{r_0}$. Hence the family
$\{z_0 z^{-1}c_\gs \colon \gs < \goo \}$ satisfies ($\dag\dag$) and
also satisfies $z_0 z^{-1}c_\gs \equiv
 z_0e_0^\gs
 \pmod{r_0}$; but this is impossible by construction.
                                                    \fin


\section*{Consequences}
Now we consider some of the general consequences, for the question of
the existence of non-standard uniserials, of the results of the
previous sections. First of all, we can construct non-standard uniserial
modules associated to any residue field of cardinality $\leq \hao$.

\begin{proposition}
\label{resfld}
(i) For any countable field $K$, there exists a valuation domain $R$
of cardinality $\hao$  with
residue field $K$ such that there is an explicitly non-standard
uniserial module of type $Q/R$.

(ii) For any field $K$ of cardinality $\leq \hao$,  there exists a valuation domain
 $R$ of cardinality $\hao$ with
residue field $K$ such that there is a non-standard
uniserial module of type $Q/R$.

\end{proposition}

\proof
Part (i) is an immediate consequence of Theorem~\ref{exist} and
Theorem~\ref{main}. Part (ii)  follows from (i)
in the case of a countable K and from the Osofsky
construction in the case of an uncountable K (cf. \cite{O1}; see also \cite[Theorem 11]{ES}).  \qed

The following improves  \cite[Corollary 15]{ES}, in that it is a theorem of
ZFC rather than a consistency result.     It shows that the condition
$\gG'(J/R) > 0$ is not necessary for the existence of a non-standard
uniserial of type $J/R$.

\begin{proposition}
\label{G'cor}
There is a valuation domain $R$ of cardinality $\hao$ such that
$\Gamma '(Q/R) = 0$ and there is a non-standard uniserial $R$-module of type
$Q/R$.
\end{proposition}

\proof
Let $R$ be such that $Q/R$ is essentially countable (cf.
Theorem~\ref{exist}).
By Theorem~\ref{gamma}(i), $\gG'(Q/R) = 0$, but
there is a non-standard uniserial $R$-module of type $Q/R$ by
Theorem~\ref{main}.
\qed

  The following sums up  some old results which
we want to combine with results proved here.

\begin{theorem}
\label{sum ES}
Suppose that $R$ is a valuation domain of cardinality $\hao$.

(i) If CH does not hold and $J/R$ is an essentially uncountable type,
 then there is a non-standard uniserial
$R$-module of type $J/R$.

(ii) If CH holds and $\Gamma(J/R) = 1$, then there is a non-standard
uniserial $R$-module of type $J/R$.

\end{theorem}

\proof
Part (i) is Theorem 8 of \cite {ES}. Part (ii) is because   the weak
diamond principle, $\Phi_{\goo}(\goo)$, is a consequence of CH (see
\cite{Sh65}) and this implies that there exists a non-standard uniserial
of type $J/R$ when $\gG(J/R) = 1$ (see \cite[Proposition 3]{ES} or
\cite{FrG}). \qed

Now we can completely
handle the cases when either CH fails, or $\Gamma = 1$.

\begin{theorem}
\label{not CH}
If CH does not hold, then for every valuation domain $R$ of cardinality
$\hao$, and every type $J/R$ such that $\Gamma(J/R) \neq 0$,
 there is a non-standard uniserial $R$-module
of type $J/R$.
\end{theorem}

\proof
Use Theorem~\ref{main} for the essentially countable case, and
Theorem~\ref{sum ES}(i) otherwise.  \qed
\noindent
{\bf Remark. }  \ This result shows that CH is needed for the independence
result  in \cite[Thm. 14]{ES}.
\begin{theorem}
\label{G=1}
For every valuation domain $R$ of cardinality $\hao$ and every type
$J/R$, if $\Gamma(J/R) = 1$, then there is a non-standard uniserial $R$-module
of type $J/R$.
\end{theorem}

\proof
If CH fails, use the previous theorem. If CH holds, use Theorem~\ref{sum
ES}(ii).  \qed

\noindent
{\bf Remark.  } \ Osofsky's original conjecture (\cite[(9), p.
164]{O1},
restricted to valuation domains of
cardinality $\hao$, said  --- in our notation --- that there is a
non-standard uniserial $R$-module of type $J/R$ if and only if $\gG(J/R)
= 1$. This is now seen to be true assuming $\neg$CH. On the other hand,
 it cannot be
true in this form assuming CH, since CH implies the weak diamond
principle for some co-stationary subsets of $\goo$ (cf.
\cite[VI.1.10]{EM}) Indeed, as in the proof
of \cite[Prop. 3]{ES} it is possible to construct $R$ with a
type $J/R$ where $\gG(J/R) =\tilde{S}$
and $\Phi_{\goo}(S)$ holds; so there is a non-standard
uniserial $R$-module of type $J/R$ (cf. \cite[proof of Prop.
3]{ES}). On the other hand, to construct such an $R$  one
has to begin with the stationary and co-stationary set $S$, so such
rings will not be ``natural'', i.e. will not be ones ordinarily met in
algebraic contexts.

Recall, from \cite{ES}, that it is in the case when the hypotheses of
the previous theorems fail --- i.e., when CH holds and $\gG(J/R) < 1$
(and non-zero) --- that the independence phenomena occur.

\section*{Absoluteness}
Finally, let us  briefly discuss absoluteness.
 Consider Theorem~\ref{not
CH}; if         CH fails and
$\gG(J/R) \neq 0$,  we always have a non-standard uniserial
module of type $J/R$,   but there are two separate constructions involved. In
one case, when $J/R$ is essentially countable, we construct an
explicitly non-standard uniserial. If the universe of set theory is
extended to a larger universe (with the same $\hao$) this module remains
non-standard because the special family of functions
 remains a special family for $U$ in the
extension of the universe.
 In the essentially
uncountable case  we use the fact that $\gG'(J/R) = 1$ ( \cite[Theorem
8]{ES}) and construct our non-standard uniserial $U$ as in \cite[Theorem
7]{ES}. In this case too $U$ remains non-standard in an extension of the
universe (preserving $\hao$). The reason here is more subtle; relative
to a fixed $\goo$-filtration of $R$ by subrings, $N_\ga$, the
$e_\gs^\gr$ we construct satisfy the following property
for every $ \gd \in \lim(\goo)$ and every $ c \in R^\ast$:

\labpar{\#_{c,\gd}}{$  \forall \langle c_\gs
\colon \gs < \gd \rangle \in {}^\gd N_\gd^\ast[
 \exists \gs <
\gd(\forall t \in R(c - c_\gs e_\gd^\gd \neq r_\gs
t))]$}

It is a theorem of ZFC that if ($\#_{c,\gd}$) holds for $U$ (defined by the
$e_\gs^\gr$) for all $c, \gd$, then $U$ is non-standard.
Now ($\#_{c,\gd}$) is, by a coding argument, a $\Pi_1^1$ statement
 (with parameters in the ground model) about $\go$.
 Hence, by a theorem of Mostowski (cf. \cite[Thm. 7.13, p. 160]{Dr}),
it remains true in an extension of the universe, so $U$ remains
non-standard.

On the other hand, in the proof of Theorem~\ref{G=1}, there is one
additional case: when  $J/R$ is
essentially uncountable  and $\gG'(J/R) = 0$ (so CH holds).
 In this case the existence of a non-standard
uniserial is proved using the weak diamond principle, which is a
consequence of CH. Here the $U$ we
construct may not remain non-standard in an extension of the universe.
Consider for example that $R$ is constructed as in \cite[Theorem
14]{ES}, but with $\gG(J/R) = 1$. If  $\oP$ is the forcing defined in the
proof there, then $\oP$ is proper, so it preserves
$\hao$ and, moreover, in the $\oP$-generic extension $U$ is standard. (Of
course, in the generic extension we can construct another non-standard
module.)


\end{document}